\documentclass[a4paper]{article}
\usepackage{pdfsync}
\usepackage[utf8]{inputenc}
\usepackage{graphicx}
\usepackage{amssymb,amsfonts,amsmath,amsthm}
\usepackage{hyperref}
\usepackage{multirow}
\usepackage{verbatim}
\usepackage{color}
\usepackage{booktabs}
\usepackage[sectionbib,round]{natbib}
\usepackage{textcomp}
\usepackage{upquote}
\usepackage{bm}
\usepackage{dirtree}

\theoremstyle{plain}

\theoremstyle{definition}
\newtheorem{definition}{Definition}[section]

\newtheorem{remark}{Remark}[section]

\begin{document}

\title{\textbf{Chance constrained directional models in stochastic data envelopment analysis}}


\author{V. J. Bol\'os$^1$, R. Ben\'{\i}tez$^1$, V. Coll-Serrano$^2$ \\ \\
{\small $^1$ Dpto. Matem\'aticas para la Econom\'{\i}a y la Empresa, Facultad de Econom\'{\i}a.} \\
{\small $^2$ Dpto. Econom\'{\i}a Aplicada, Facultad de Econom\'{\i}a.} \\
{\small Universidad de Valencia. Avda. Tarongers s/n, 46022 Valencia, Spain.} \\
{\small e-mail\textup{: \texttt{vicente.bolos@uv.es}, \texttt{rabesua@uv.es}, \texttt{vicente.coll@uv.es}}} \\}

\date{February 2024}

\maketitle

\begin{abstract}
We construct a new family of chance constrained directional models in stochastic data envelopment analysis, generalizing the deterministic directional models and the chance constrained radial models. We prove that chance constrained directional models define the same concept of stochastic efficiency as the one given by chance constrained radial models and, as a particular case, we obtain a stochastic version of the generalized Farrell measure. Finally, we give some examples of application of chance constrained directional models with stochastic and deterministic directions, showing that inefficiency scores obtained with stochastic directions are less or equal than those obtained considering deterministic directions whose values are the means of the stochastic ones.
\end{abstract}

\section{Introduction}

Data envelopment analysis (DEA) \citep{Charnes1978} is a non-parametric technique used to measure the relative efficiency of a homogeneous set of decision-making units (DMUs) that use multiple inputs to obtain multiple outputs. Using mathematical programming methods, DEA allows the identification of the best practice frontier of the production possibility set, determined by DMUs which are qualified as efficient. On the other hand, DMUs that move away from this frontier are considered inefficient.

Radial models were introduced by \citet{Charnes1978, Charnes1979, Charnes1981} for constant returns to scale (CCR models), and \citet{Banker1984} for variable returns to scale (BCC models). In general, they can be either input or output-oriented. In the former case, we look for determining the maximal proportionate reduction of inputs allowed by the production possibility set, while maintaining the current output level. On the other hand, in the output-oriented case, we want to find the maximal proportionate increase of outputs while keeping the current input consumption.

One of the drawbacks of radial models is that these reductions (of inputs) or increases (of outputs) have to be made in the same proportion. Therefore, \citet{Fare1978} defined non-radial models that allow non-proportional movements in inputs and outputs. Subsequently, \citet{Charnes1985} introduced the so-called additive models, which are also non-radial and do not distinguish between orientations. Later, \citet{Pastor1999} and \citet{Tone2001} constructed the slacks-based measure of efficiency (SBM) models, with the same characteristics as additive models but providing an efficiency score. However, all these models compute efficiency with respect to the target located at the ``farthest'' point within the portion of the efficient frontier that dominates the evaluated DMU and, therefore, may be inappropriate in the case of wanting to calculate a ``close'' efficient target. Some authors adapted these non-radial models to the search of the closest targets. For example, \citet{Apa2007} adapted several models and constructed the mADD model as the minimum-distance version of the additive model. In turn, \citet{Tone2010} developed the SBM-max model. Nevertheless, these models are not weakly monotonic in general, in the sense that a DMU can get a worse efficiency score than another DMU it outperforms. 

Other models allowing non-proportional movements in inputs and outputs were introduced by \citet{Luenberger1992}. Later, \citet{Chambers1996,Chambers1998} and \citet{Briec1997} introduced the directional models. These models, apart from being weakly monotonic, generalize the classical radial models and can measure efficiency in the full input-output space. The main goal of directional models is that the user can change the proportion in which the inputs and outputs are modified, thus defining a custom orientation taking account of particularities of the market and characterizing the criteria of management chosen by the producer.

Nevertheless, all these models assume a deterministic behaviour of the inputs and outputs of the DMUs, and are not appropriate in a stochastic framework. In this sense, stochastic DEA was developed in two main directions (see \citet{Olesen2016} for a complete review on stochastic DEA). On one hand, \citet{Banker1993} initiated one approach that included statistical axioms defining a statistical model and a sampling process into the DEA framework. On the other hand, chance constrained DEA \citep{Lan1993, Olesen1995} allows for uncertainty in the inputs and outputs, ensuring that the efficiency scores are obtained with a given degree of confidence and hence, providing a more robust and realistic assessment of efficiency compared to deterministic DEA. \citet{Cooper1996,Cooper1998,Cooper2002} developed the chance constrained version of some radial models, but, up to our knowledge, there are no chance constrained directional models in the literature. Therefore, in the same way that directional models generalize the classical radial models in the deterministic DEA, the objective of this work is to define chance constrained directional models, generalizing the existing chance constrained radial models.

Finally, it is worth mentioning that fuzzy DEA is another recent approach to DEA analysis presenting uncertainty in inputs and outputs, with Kao-Liu \citep{Kao2000}, Guo-Tanaka \citep{Guo2001} and possibilistic \citep{Leon2003} models as the most relevant. More recently, some chance constrained fuzzy DEA models has been developed by \citet{Tavana2013}. However, fuzzy DEA could be considered too limited in a stochastic framework because it does not take into account the probabilistic distribution of the data, resulting in a potential loss of significant details.

This paper is organized as follows. In Section \ref{sec:dir}, we review directional models and adapt the notation to be able to define, in Section \ref{sec:ccd}, their chance constrained versions in which inputs and outputs become stochastic. We have developed two types of chance constrained directional models: with \emph{stochastic directions} and with \emph{deterministic directions}. Models with stochastic directions are introduced in Section \ref{sec:ccdsto} and they generalize the existing chance constrained radial models. Moreover, they serve to define a stochastic version of the generalized Farrell measure. On the other hand, models with deterministic directions are defined analogously in Section \ref{sec:ccddet}. We give some examples of application of these models in Section \ref{sec:ex} and finally, we present some conclusions and raise some open issues in Section \ref{sec:con}.

In general, we denote vectors by bold-face letters and they are considered as column vectors unless otherwise stated. The elements of a vector are denoted by the same letter as the vector, but unbolded and with subscripts. The $0$-vector is denoted by $\bm{0}$ and the context determines its dimension.

\section{Directional models}
\label{sec:dir}

We consider $\mathcal{D}=\left\{ \textrm{DMU}_1, \ldots ,\textrm{DMU}_n \right\} $ a set of $n$ DMUs with $m$ inputs and $s$ outputs. Matrices $X=(x_{ij})$ and $Y=(y_{rj})$ are the \emph{input} and \emph{output data matrices}, respectively, where $x_{ij}$ and $y_{rj}$ denote the $i$th input and $r$th output of the $j$th DMU. We also assume that $x_{ij}$ and $y_{rj}$ are all strictly positive. Given $\textrm{DMU}_o\in \mathcal{D}$, we define the column vectors $\mathbf{x}_o=(x_{1o},\ldots,x_{mo})^\top$ and $\mathbf{y}_o=(y_{1o},\ldots,y_{so})^\top$.
Although we are going to consider constant returns to scale (CRS) in all programs, different returns to scale can be easily considered by adding the corresponding constraints: $\mathbf{e}\bm{\lambda}=1$ for variable returns to scale (VRS), $0\leq \mathbf{e}\bm{\lambda}\leq 1$ for non-increasing returns to scale (NIRS), $\mathbf{e}\bm{\lambda}\geq 1$ for non-decreasing returns to scale (NDRS) or $L\leq \mathbf{e}\bm{\lambda}\leq U$ for generalized returns to scale (GRS), with $0\leq L\leq 1$ and $U\geq 1$, where $\bm{\lambda}\in \mathbb{R}^n$ and $\mathbf{e}=(1,\ldots,1)\in \mathbb{R}^n$ is the all-ones row vector.

Directional models were introduced in \cite{Chambers1996,Chambers1998}. Given $\textrm{DMU}_o\in \mathcal{D}$, the associated linear program under CRS and its second stage are given by
\begin{equation}
\label{eq:dir}
\def\arraystretch{1.2}
\begin{array}{ll}
\begin{array}[t]{rl}
\textrm{(a) }\max \limits_{\beta, \bm{\lambda}} & \beta \\
\text{s.t.} & \beta \mathbf{g}^-+X\bm{\lambda}\leq \mathbf{x}_o, \\
& -\beta \mathbf{g}^++Y\bm{\lambda} \geq \mathbf{y}_o, \\
& \bm{\lambda}\geq \mathbf{0},
\end{array}
&\quad
\begin{array}[t]{rl}
\textrm{(b) }\max \limits_{\bm{\lambda}, \mathbf{s}^-,\mathbf{s}^+} & \mathbf{w}^-\mathbf{s}^-+\mathbf{w}^+\mathbf{s}^+ \\
\text{s.t.} & X\bm{\lambda}+\mathbf{s}^-= \mathbf{x}_o-\beta ^*\mathbf{g}^-, \\
& Y\bm{\lambda} -\mathbf{s}^+= \mathbf{y}_o+\beta ^*\mathbf{g}^+, \\
& \bm{\lambda}\geq \mathbf{0},\,\, \mathbf{s}^-\geq \mathbf{0},\,\, \mathbf{s}^+\geq \mathbf{0},
\end{array}
\end{array}
\end{equation}
where $\bm{\lambda}\in \mathbb{R}^n$ and the \emph{slacks} $\mathbf{s}^-\in \mathbb{R}^m$ and $\mathbf{s}^+\in \mathbb{R}^s$ are non-negative column vectors, $\mathbf{g}=(-\mathbf{g}^-,\mathbf{g}^+)\neq \mathbf{0}$ is a preassigned \emph{direction} (with $\mathbf{g}^-\in \mathbb{R}^m$ and $\mathbf{g}^+\in\mathbb{R}^s$ non-negative column vectors), while the \emph{weights} $\mathbf{w}^-\in \mathbb{R}^m$ and $\mathbf{w}^+\in \mathbb{R}^s$ are strictly positive row vectors. Moreover, in \eqref{eq:dir} (b), $\beta ^*$ is the optimal objective value of the first stage program \eqref{eq:dir} (a), which is always greater than or equal to $0$.

The advantage of directional models is that the way DMU$_o$ is projected onto the efficient frontier can be customized by using a direction $\mathbf{g}$ according to the criteria of management chosen by the producer \citep{Briec1997}.

\begin{definition}
\label{def:1}
$\textrm{DMU}_o\in \mathcal{D}$ is \emph{efficient} if and only if the optimal objectives $\beta ^*$ and $\mathbf{w}^-\mathbf{s}^{-*}+\mathbf{w}^+\mathbf{s}^{+*}$ in \eqref{eq:dir} are $0$, where $\mathbf{s}^{-*}, \mathbf{s}^{+*}$ are optimal slacks. Moreover, we say that $\textrm{DMU}_o$ is \emph{weakly efficient} if $\beta ^*=0$ and it is not efficient.
\end{definition}

\begin{remark}
\label{teo:1}
The concept of efficiency given in Definition \ref{def:1} does not depend on the direction $\mathbf{g}$ because $\beta ^*=0$ and hence, according to \eqref{eq:dir}, the optimal solution remains optimal if we change the direction. So, if $\textrm{DMU}_o$ is efficient for a given direction, then it is efficient for any direction.
\end{remark}

\begin{remark} We consider that the weights $\mathbf{w}^-$ and $\mathbf{w}^+$ are strictly positive and hence, $\mathbf{w}^-\mathbf{s}^{-*}+\mathbf{w}^+\mathbf{s}^{+*}=0$ if and only if $\mathbf{s}^{-*}=\mathbf{0}$ and $\mathbf{s}^{+*}=\mathbf{0}$. However, allowing zero weights is useful in some cases. For example, we can introduce \emph{non-discretionary} inputs and outputs (which are exogenously fixed, see \citet[Chapter 7]{Cooper2007}) by setting the corresponding directions and weights $g_i^-=w_i^-=0$ (for the $i$th input) or $g_r^+=w_r^+=0$ (for the $r$th output). In this case, the optimal slacks of non-discretionary variables are not taken into account in the efficiency evaluation.
\end{remark}

Nevertheless, in order to construct the chance constrained version of the model in Section \ref{sec:ccd}, it is convenient to write the second stage program as
\begin{equation}
\label{eq:dir2}
\begin{array}[t]{rl}
\max \limits_{\bm{\lambda}, \mathbf{s}^-,\mathbf{s}^+} & \mathbf{w}^-\mathbf{s}^-+\mathbf{w}^+\mathbf{s}^+ \\
\text{s.t.} & X\bm{\lambda}+\mathbf{s}^-\leq \mathbf{x}_o-\beta ^*\mathbf{g}^-, \\
& Y\bm{\lambda} -\mathbf{s}^+\geq \mathbf{y}_o+\beta ^*\mathbf{g}^+, \\
& \bm{\lambda}\geq \mathbf{0},\,\, \mathbf{s}^-\geq \mathbf{0},\,\, \mathbf{s}^+\geq \mathbf{0}.
\end{array}
\end{equation}
Obviously, programs \eqref{eq:dir} (b) and \eqref{eq:dir2} are equivalent, since any optimal solution of \eqref{eq:dir2} should also satisfy the constraints of \eqref{eq:dir} (b), and vice versa.

The one-stage version of \eqref{eq:dir} (a) and \eqref{eq:dir2} is given by
\begin{equation}
\label{eq:dir1}
\begin{array}[t]{rl}
\max \limits_{\beta ,\bm{\lambda}, \mathbf{s}^-,\mathbf{s}^+} & \beta +\varepsilon \left( \mathbf{w}^-\mathbf{s}^-+\mathbf{w}^+\mathbf{s}^+\right) \\
\text{s.t.} & \beta \mathbf{g}^-+X\bm{\lambda}+\mathbf{s}^-\leq \mathbf{x}_o, \\
& -\beta \mathbf{g}^++Y\bm{\lambda} -\mathbf{s}^+\geq \mathbf{y}_o, \\
& \bm{\lambda}\geq \mathbf{0},\,\, \mathbf{s}^-\geq \mathbf{0},\,\, \mathbf{s}^+\geq \mathbf{0},
\end{array}
\end{equation}
where $\varepsilon$ a positive non-Archimedean infinitesimal. Although the use of non-Archimedean infinitesimals is not numerically rigorous, in the case of DEA models with two stages is very useful because it synthesizes the two stages in a single program and you can always do the reverse process to recover the programs of each stage.

Let us consider the directions $\mathbf{g}^-=D^-\mathbf{x}_o$ and $\mathbf{g}^+=D^+\mathbf{y}_o$, where the matrices $D^-=\mathrm{diag}(d^-_1,\ldots ,d^-_m)$ and $D^+=\mathrm{diag}(d^+_1,\ldots ,d^+_s)$ are diagonal with $d^-_1,\ldots ,d^-_m,d^+_1,\ldots ,d^+_s\geq 0$. Then \eqref{eq:dir1} becomes
\begin{equation}
\label{eq:dirD1}
\def\arraystretch{1.2}
\begin{array}[t]{rl}
\max \limits_{\beta ,\bm{\lambda}, \mathbf{s}^-,\mathbf{s}^+} & \beta +\varepsilon \left( \mathbf{w}^-\mathbf{s}^-+\mathbf{w}^+\mathbf{s}^+\right) \\
\text{s.t.} & \Theta ^-(\beta )\mathbf{x}_o-X\bm{\lambda}-\mathbf{s}^-\geq \mathbf{0}, \\
& \Theta ^+(\beta )\mathbf{y}_o-Y\bm{\lambda}+\mathbf{s}^+\leq \mathbf{0}, \\
& \bm{\lambda}\geq \mathbf{0},\,\, \mathbf{s}^-\geq \mathbf{0},\,\, \mathbf{s}^+\geq \mathbf{0},
\end{array}
\end{equation}
where 
\begin{equation}
\label{eq:Theta}
\Theta ^-(\beta )=I_m-\beta D^-,\qquad \Theta ^+(\beta )=I_s+\beta D^+,
\end{equation}
with $I_m$, $I_s$ the $m\times m$ and $s\times s$ identity matrices, respectively. The optimal value $\beta ^*$ given by \eqref{eq:dirD1} coincides with the \emph{oriented Farrell proportional distance} defined by \cite{Briec1997}, with \emph{orientation} $\left( \mathbf{d}^-,\mathbf{d}^+\right) $, where $\mathbf{d}^-=(d^-_1,\ldots ,d^-_m)$ and $\mathbf{d}^+=(d^+_1,\ldots ,d^+_s)$. Hence, by Definition \ref{def:1}, $\textrm{DMU}_o$ is efficient if and only if the optimal objective value of \eqref{eq:dirD1} is $0$, regardless of the orientation. 

As a particular case, if $D^-=I_m$ (i.e. $d^-_1,\ldots ,d^-_m=1$), and $D^+=0$ (i.e. $d^+_1,\ldots ,d^+_s=0$), then we can write \eqref{eq:dirD1} as
\begin{equation}
\label{eq:rad1a}
\def\arraystretch{1.2}
\begin{array}{rl}
\min \limits_{\theta ,\bm{\lambda}, \mathbf{s}^-,\mathbf{s}^+} & \theta -\varepsilon \left( \mathbf{w}^-\mathbf{s}^-+\mathbf{w}^+\mathbf{s}^+\right) \\
\text{s.t.} & \theta \mathbf{x}_o-X\bm{\lambda}-\mathbf{s}^-\geq \mathbf{0}, \\
& Y\bm{\lambda}-\mathbf{s}^+\geq \mathbf{y}_o, \\
& \bm{\lambda}\geq \mathbf{0},\,\, \mathbf{s}^-\geq \mathbf{0},\,\, \mathbf{s}^+\geq \mathbf{0},
\end{array}
\end{equation}
where $\theta =1-\beta $. Program \eqref{eq:rad1a} corresponds to the input-oriented CRS radial model, also known as CCR model. So, according to Remark \ref{teo:1} and considering strictly positive weights, the concept of efficiency given in Definition \ref{def:1} coincides with the classical \textit{CCR-efficiency} \citep[Chapter 3]{Cooper2007}.

On the other hand, if $D^-=0$ (i.e. $d^-_1,\ldots ,d^-_m=0$), and $D^+=I_s$ (i.e. $d^+_1,\ldots ,d^+_s=1$), then  we can write \eqref{eq:dirD1} as
\begin{equation}
\label{eq:rad1b}
\def\arraystretch{1.2}
\begin{array}{rl}
\max \limits_{\phi ,\bm{\lambda}, \mathbf{s}^-,\mathbf{s}^+} & \phi +\varepsilon \left( \mathbf{w}^-\mathbf{s}^-+\mathbf{w}^+\mathbf{s}^+\right) \\
\text{s.t.} & X\bm{\lambda}+\mathbf{s}^-\leq \mathbf{x}_o, \\
& \phi \mathbf{y}_o-Y\bm{\lambda}+\mathbf{s}^+\leq \mathbf{0}, \\
& \bm{\lambda}\geq \mathbf{0},\,\, \mathbf{s}^-\geq \mathbf{0},\,\, \mathbf{s}^+\geq \mathbf{0},
\end{array}
\end{equation}
where $\phi =1+\beta $. Program \eqref{eq:rad1b} is the output-oriented version of the CCR model.

Moreover, if $D^-=I_m$ and $D^+=I_s$ (i.e. $d^-_1,\ldots ,d^-_m,d^+_1,\ldots ,d^+_s=1$), then model \eqref{eq:dirD1} is non-oriented and $\beta ^*$ coincides with the \textit{generalized Farrell measure} \citep{Briec1997}.

\section{Chance constrained directional models}
\label{sec:ccd}

We consider $\mathcal{D}=\left\{ \textrm{DMU}_1, \ldots ,\textrm{DMU}_n \right\} $ a set of $n$ DMUs with $m$ stochastic inputs and $s$ stochastic outputs. Matrices $\tilde{X}=(\tilde{x}_{ij})$ and $\tilde{Y}=(\tilde{y}_{rj})$ are the \emph{input} and \emph{output data matrices}, respectively, where $\tilde{x}_{ij}$ and $\tilde{y}_{rj}$ represent the $i$th input and $r$th output of the $j$th DMU. Moreover, we denote by $X=(x_{ij})$ and $Y=(y_{rj})$ their expected values, which we assume to be strictly positive. Given $\textrm{DMU}_o\in \mathcal{D}$, we define the column vectors $\tilde{\mathbf{x}}_o=(\tilde{x}_{1o},\ldots ,\tilde{x}_{mo})^\top$, $\tilde{\mathbf{y}}_o=(\tilde{y}_{1o},\ldots ,\tilde{y}_{so})^\top$, $\mathbf{x}_o=(x_{1o},\ldots ,x_{mo})^\top$ and $\mathbf{y}_o=(y_{1o},\ldots ,y_{so})^\top$.
Although we suppose CRS in all programs, different returns to scale can be easily considered by adding the corresponding constraints: $\mathbf{e}\bm{\lambda}=1$ (VRS), $0\leq \mathbf{e}\bm{\lambda}\leq 1$ (NIRS), $\mathbf{e}\bm{\lambda}\geq 1$ (NDRS) or $L\leq \mathbf{e}\bm{\lambda}\leq U$ (GRS), with $0\leq L\leq 1$ and $U\geq 1$.

In general, there are two methods for constructing a chance constrained model from a deterministic model in which inputs and outputs become stochastic, according to \citet[Chapter 4]{Amir2022}. On the one hand, a \emph{P-model} is based on the highest probability of occurrence of the objective function. On the other hand, an \emph{E-model} is based on the mathematical expectation to obtain the expected value of the objective function. In both cases, constraints must be satisfied with a probability greater than a certain threshold.

In particular, we want to adapt the deterministic directional models \eqref{eq:dir1} or \eqref{eq:dirD1} to their chance constrained versions in which inputs and outputs become stochastic. Since the objective function remains deterministic, we are going to apply the E-model approach.

Note that, although programs \eqref{eq:dir1} and \eqref{eq:dirD1} are equivalent, there is a relevant difference if we want to adapt them to the stochastic case: the directions $\mathbf{g}^-$ and $\mathbf{g}^+$ from \eqref{eq:dir1} remain deterministic but, on the other hand, the directions $\mathbf{g}^-=D^-\tilde{\mathbf{x}}_o$ and $\mathbf{g}^+=D^+\tilde{\mathbf{y}}_o$ from \eqref{eq:dirD1} become stochastic because they depend on the (stochastic) inputs and outputs. So, we are mostly interested in adapting program \eqref{eq:dirD1} because it generalizes the radial models, in which the radial direction depends on the inputs or outputs. Hence, we will first adapt \eqref{eq:dirD1} in Section \ref{sec:ccdsto} and later, adapt \eqref{eq:dir1} in Section \ref{sec:ccddet}.

\subsection{Stochastic directions}
\label{sec:ccdsto}

In this section, we are going to apply the E-model approach to program \eqref{eq:dirD1}, considering that inputs and outputs become stochastic. Hence, given $0<\alpha <1$, the corresponding chance constrained E-model can be written as
\begin{equation}
\label{eq:ccdirD21a}
\def\arraystretch{1.8}
\begin{array}[t]{rl}
\max \limits_{\beta, \bm{\lambda},\mathbf{s}^-,\mathbf{s}^+} & \beta +\varepsilon \left( \mathbf{w}^-\mathbf{s}^-+\mathbf{w}^+\mathbf{s}^+\right) \\
\text{s.t.} & P\left\{ \left( \Theta ^-(\beta )\tilde{\mathbf{x}}_o-\tilde{X}\bm{\lambda}-\mathbf{s}^-\right) _i\geq 0\right\} \geq 1-\alpha,\quad i=1,\ldots ,m, \\
& P\left\{ \left( \Theta ^+(\beta )\tilde{\mathbf{y}}_{o}-\tilde{Y}\bm{\lambda}+\mathbf{s}^+\right) _r\leq 0\right\} \geq 1-\alpha,\quad r=1,\ldots ,s, \\
& \bm{\lambda}\geq \mathbf{0},\,\, \mathbf{s}^-\geq \mathbf{0},\,\, \mathbf{s}^+\geq \mathbf{0},
\end{array}
\end{equation}
where $P$ denotes the probability function and $\Theta ^-(\beta )$, $\Theta ^+(\beta )$ are given by \eqref{eq:Theta}.

\begin{definition}
\label{def:2}
$\textrm{DMU}_o\in \mathcal{D}$ is \emph{$\alpha$-stochastically efficient} if and only if the optimal objective value of \eqref{eq:ccdirD21a} is $0$, i.e. $\beta ^*=0$ and  $\mathbf{w}^-\mathbf{s}^{-*}+\mathbf{w}^+\mathbf{s}^{+*}=0$. Moreover, we say that $\textrm{DMU}_o$ is \emph{$\alpha$-stochastically weakly efficient} if $\beta ^*=0$ and it is not $\alpha$-stochastically efficient.
\end{definition}

\begin{remark}
\label{teo:2}
The concept of stochastic efficiency given in Definition \ref{def:2} does not depend on the orientation $(\mathbf{d}^-,\mathbf{d}^+)$ because $\beta ^*=0$ and hence, by \eqref{eq:Theta}, $\Theta^-(\beta ^*)=I_m$ and $\Theta^+(\beta ^*)=I_s$ do not depend on the orientation. So, according to \eqref{eq:ccdirD21a}, optimal solutions remain optimal if we change the orientation, concluding that  if $\textrm{DMU}_o$ is $\alpha$-stochastically efficient for a given orientation, then it is $\alpha$-stochastically efficient for any orientation.
\end{remark}

For constraints in \eqref{eq:ccdirD21a}, we have
\begin{align}
\label{eq:paso1i}
P\left\{ \left( \Theta ^-(\beta )\tilde{\mathbf{x}}_o-\tilde{X}\bm{\lambda}\right) _i\leq s_i^-\right\} &\leq \alpha ,\quad i=1,\ldots ,m,\\
\label{eq:paso1o} P\left\{ \left( \tilde{Y}\bm{\lambda}-\Theta ^+(\beta )\tilde{\mathbf{y}}_o\right) _r\leq s_r^+\right\} &\leq \alpha ,\quad r=1,\ldots ,s.
\end{align}
Let us assume that inputs and outputs are correlated random variables following multivariate normal distributions with means $E(\tilde{x}_{ij})=x_{ij}$ and $E(\tilde{y}_{rj})=y_{rj}$. The use of the normal distribution is due to its great versatility and it is discussed by \citet{Cooper1996}. Then, we can define standard normal random variables by
\begin{equation*}
\tilde{Z}^-_i=\dfrac{\left( \Theta ^-(\beta )\tilde{\mathbf{x}}_o-\tilde{X}\bm{\lambda}\right) _i-\Big( \Theta ^-(\beta )\mathbf{x}_o-X\bm{\lambda}\Big) _i}{\sigma ^-_i\left( \beta, \bm{\lambda}\right)},\qquad i=1,\ldots ,m,
\end{equation*}
\begin{equation*}
\tilde{Z}^+_r=\dfrac{\left( \tilde{Y}\bm{\lambda}-\Theta ^+(\beta )\tilde{\mathbf{y}}_o\right) _r-\Big( Y\bm{\lambda}-\Theta ^+(\beta )\mathbf{y}_o\Big) _r}{\sigma ^+_r\left( \beta, \bm{\lambda}\right)},\qquad r=1,\ldots ,s,
\end{equation*}
where
\begin{align}
\nonumber
\left( \sigma ^-_i\left( \beta, \bm{\lambda}\right) \right) ^2 =& \,\mathrm{Var}\left( \Theta ^-(\beta )\tilde{\mathbf{x}}_o-\tilde{X}\bm{\lambda}\right) _i=\mathrm{Var}\left( (1-\beta d^-_i)\tilde{x}_{io}-\sum _{j=1}^n \lambda _j\tilde{x}_{ij}\right) \\
\label{eq:sigma1}
=& \sum _{j,q=1}^n\lambda _j\lambda _q\mathrm{Cov}(\tilde{x}_{ij},\tilde{x}_{iq})-2(1-\beta d^-_i)\sum _{j=1}^n\lambda _j\mathrm{Cov}(\tilde{x}_{ij},\tilde{x}_{io}) \\
\nonumber & +(1-\beta d^-_i)^2\mathrm{Var}(\tilde{x}_{io}),\quad i=1,\ldots ,m,\\
\nonumber \left( \sigma ^+_r\left( \beta, \bm{\lambda}\right)\right) ^2 =& \,\mathrm{Var}\left( \tilde{Y}\bm{\lambda}-\Theta ^+(\beta )\tilde{\mathbf{y}}_o\right) _r=\mathrm{Var}\left( \sum _{j=1}^n \lambda _j\tilde{y}_{rj}-(1+\beta d^+_r)\tilde{y}_{ro}\right) \\
\label{eq:sigma2}
=& \sum _{j,q=1}^n\lambda _j\lambda _q\mathrm{Cov}(\tilde{y}_{rj},\tilde{y}_{rq})-2(1+\beta d^+_r)\sum _{j=1}^n\lambda _j\mathrm{Cov}(\tilde{y}_{rj},\tilde{y}_{ro}) \\
\nonumber & +(1+\beta d^+_r)^2\mathrm{Var}(\tilde{y}_{ro}),\quad r=1,\ldots ,s.
\end{align}
Hence, \eqref{eq:paso1i} and \eqref{eq:paso1o} become
\begin{align}
\label{eq:paso2i}
P\left\{ \tilde{Z}^-_i \leq \dfrac{s_i^--\left( \Theta ^-(\beta )\mathbf{x}_o-X\bm{\lambda}\right) _i}{\sigma ^-_i\left( \beta, \bm{\lambda}\right)}\right\} &\leq \alpha ,\quad i=1,\ldots ,m,\\
\label{eq:paso2o} P\left\{ \tilde{Z}^+_r\leq \dfrac{s_r^+-\left( Y\bm{\lambda}-\Theta ^+(\beta )\mathbf{y}_o\right) _r}{\sigma ^+_r\left( \beta, \bm{\lambda}\right)}\right\} &\leq \alpha ,\quad r=1,\ldots ,s.
\end{align}
Since $\tilde{Z}^-_i$ and $\tilde{Z}^+_r$ follow the standard normal distribution $\Phi$, we can write \eqref{eq:paso2i} and \eqref{eq:paso2o} as
\begin{align}
\label{eq:paso3i}
\dfrac{s_i^--\left( \Theta ^-(\beta )\mathbf{x}_o-X\bm{\lambda}\right) _i}{\sigma ^-_i\left( \beta, \bm{\lambda}\right)} &\leq \Phi ^{-1}(\alpha ),\quad i=1,\ldots ,m,\\
\label{eq:paso3o} \dfrac{s_r^+-\left( Y\bm{\lambda}-\Theta ^+(\beta )\mathbf{y}_o\right) _r}{\sigma ^+_r\left( \beta, \bm{\lambda}\right)} &\leq \Phi ^{-1}(\alpha ),\quad r=1,\ldots ,s.
\end{align}
Reordering and vectorizing expressions \eqref{eq:paso3i} and \eqref{eq:paso3o}, we have that the deterministic equivalent to \eqref{eq:ccdirD21a} is given by
\begin{equation}
\label{eq:ccdirD2d0}
\def\arraystretch{1.2}
\begin{array}[t]{rl}
\max \limits_{\beta ,\bm{\lambda},\mathbf{s}^-,\mathbf{s}^+} & \beta +\varepsilon \left( \mathbf{w}^-\mathbf{s}^-+\mathbf{w}^+\mathbf{s}^+\right) \\
\text{s.t.} & \Theta ^-(\beta )\mathbf{x}_o-X\bm{\lambda}-\mathbf{s}^-+\Phi ^{-1}(\alpha)\bm{\sigma} ^-(\beta ,\bm{\lambda}) \geq \mathbf{0},\\
& \Theta ^+(\beta )\mathbf{y}_{o}-Y\bm{\lambda}+\mathbf{s}^+-\Phi ^{-1}(\alpha)\bm{\sigma} ^+(\beta ,\bm{\lambda}) \leq \mathbf{0},\\
& \bm{\lambda}\geq \mathbf{0},\,\, \mathbf{s}^-\geq \mathbf{0},\,\, \mathbf{s}^+\geq \mathbf{0},
\end{array}
\end{equation}
that is equivalent to
\begin{equation}
\label{eq:ccdirD2d1}
\def\arraystretch{1.2}
\begin{array}[t]{rl}
\max \limits_{\beta ,\bm{\lambda},\mathbf{s}^-,\mathbf{s}^+} & \beta +\varepsilon \left( \mathbf{w}^-\mathbf{s}^-+\mathbf{w}^+\mathbf{s}^+\right) \\
\text{s.t.} & \Theta ^-(\beta )\mathbf{x}_o-X\bm{\lambda}-\mathbf{s}^-+\Phi ^{-1}(\alpha)\bm{\sigma} ^-(\beta ,\bm{\lambda}) = \mathbf{0},\\
& \Theta ^+(\beta )\mathbf{y}_{o}-Y\bm{\lambda}+\mathbf{s}^+-\Phi ^{-1}(\alpha)\bm{\sigma} ^+(\beta ,\bm{\lambda}) = \mathbf{0},\\
& \bm{\lambda}\geq \mathbf{0},\,\, \mathbf{s}^-\geq \mathbf{0},\,\, \mathbf{s}^+\geq \mathbf{0},
\end{array}
\end{equation}
because any optimal solution of \eqref{eq:ccdirD2d0} should also satisfy the constraints of \eqref{eq:ccdirD2d1}, and vice versa. Hence, by Definition \ref{def:2}, $\textrm{DMU}_o$ is $\alpha$-stochastically efficient if and only if the optimal objective value of \eqref{eq:ccdirD2d1} is $0$. Note that, considering $\bm{\sigma}^-,\bm{\sigma}^+$ as variables and adding \eqref{eq:sigma1}, \eqref{eq:sigma2} as constraints, program \eqref{eq:ccdirD2d1} is quadratic.

\begin{remark}
\label{rem:multi}
The chance constrained methodology ensures that the constraints are met with a certain probability, but considering each constraint separately. A more appropriate study would be to require that all the constraints are met at once with a certain probability. For example, program \eqref{eq:ccdirD21a} would look like this:
\begin{equation}
\label{eq:ccdirD21aconc}
\def\arraystretch{1.8}
\begin{array}[t]{rl}
\max \limits_{\beta, \bm{\lambda},\mathbf{s}^-,\mathbf{s}^+} & \beta +\varepsilon \left( \mathbf{w}^-\mathbf{s}^-+\mathbf{w}^+\mathbf{s}^+\right) \\
\text{s.t.} & P\left\{ 
\begin{array}{l}
\left( \Theta ^-(\beta )\tilde{\mathbf{x}}_o-\tilde{X}\bm{\lambda}-\mathbf{s}^-\geq \mathbf{0}\right) \wedge \\
\left( \Theta ^+(\beta )\tilde{\mathbf{y}}_{o}-\tilde{Y}\bm{\lambda}+\mathbf{s}^+\leq \mathbf{0}\right) 
\end{array}
\right\} \geq 1-\alpha, \\
& \bm{\lambda}\geq \mathbf{0},\,\, \mathbf{s}^-\geq \mathbf{0},\,\, \mathbf{s}^+\geq \mathbf{0},
\end{array}
\end{equation}
where $\wedge$ denotes the logical conjunction.
But, in order to find the deterministic equivalent program of \eqref{eq:ccdirD21aconc}, we should use multivariate probability distributions, for which the inverse distribution functions are not uniquely defined. An approach to overcome this problem is to consider independent variables, so that the joint probability is the product of the separate probabilities \citep{Shiraz2021}, although this method reduces the applicability of the model.
\end{remark}

As a particular case of \eqref{eq:ccdirD2d1}, if $D^-=I_m$ (i.e. $d^-_1,\ldots ,d^-_m=1$), and $D^+=0$ (i.e. $d^+_1,\ldots ,d^+_s=0$), then we can write
\begin{equation}
\label{eq:ccr2d1i}
\def\arraystretch{1.2}
\begin{array}[t]{rl}
\min \limits_{\theta ,\bm{\lambda},\mathbf{s}^-,\mathbf{s}^+} & \theta -\varepsilon \left( \mathbf{w}^-\mathbf{s}^-+\mathbf{w}^+\mathbf{s}^+\right) \\
\text{s.t.} & \theta \mathbf{x}_o-X\bm{\lambda}-\mathbf{s}^-+\Phi ^{-1}(\alpha)\bm{\sigma} ^-(\theta ,\bm{\lambda}) =\mathbf{0},\\
& Y\bm{\lambda}-\mathbf{s}^++\Phi ^{-1}(\alpha)\bm{\sigma} ^+(\bm{\lambda}) =\mathbf{y}_{o},\\
& \bm{\lambda}\geq \mathbf{0},\,\, \mathbf{s}^-\geq \mathbf{0},\,\, \mathbf{s}^+\geq \mathbf{0},
\end{array}
\end{equation}
where $\theta =1-\beta $, and
\begin{align*}
\left( \sigma ^-_i\left( \theta, \bm{\lambda}\right)\right) ^2 =& \sum _{j,q=1}^n\lambda _j\lambda _q\mathrm{Cov}(\tilde{x}_{ij},\tilde{x}_{iq})-2\theta\sum _{j=1}^n\lambda _j\mathrm{Cov}(\tilde{x}_{ij},\tilde{x}_{io})+\theta ^2\mathrm{Var}(\tilde{x}_{io}),\\
& i=1,\ldots ,m,\\
\left( \sigma ^+_r\left( \bm{\lambda}\right)\right) ^2 =& \sum _{j,q=1}^n\lambda _j\lambda _q\mathrm{Cov}(\tilde{y}_{rj},\tilde{y}_{rq})-2\sum _{j=1}^n\lambda _j\mathrm{Cov}(\tilde{y}_{rj},\tilde{y}_{ro})+\mathrm{Var}(\tilde{y}_{ro}),\\
& r=1,\ldots ,s.
\end{align*}
Program \eqref{eq:ccr2d1i} is the chance constrained E-model corresponding to the input-oriented CCR model \eqref{eq:rad1a}, firstly introduced with some simplifications by \citet{Lan1993}.

On the other hand, if $D^-=0$ (i.e. $d^-_1,\ldots ,d^-_m=0$), and $D^+=I_s$ (i.e. $d^+_1,\ldots ,d^+_s=1$), then we can write \eqref{eq:ccdirD2d1} as
\begin{equation}
\label{eq:ccr2d1o}
\def\arraystretch{1.2}
\begin{array}[t]{rl}
\max \limits_{\phi ,\bm{\lambda},\mathbf{s}^-,\mathbf{s}^+} & \phi +\varepsilon \left( \mathbf{w}^-\mathbf{s}^-+\mathbf{w}^+\mathbf{s}^+\right) \\
\text{s.t.} & X\bm{\lambda}+\mathbf{s}^--\Phi ^{-1}(\alpha)\bm{\sigma} ^-(\bm{\lambda}) =\mathbf{x}_o,\\
& \phi \mathbf{y}_{o}-Y\bm{\lambda}+\mathbf{s}^+-\Phi ^{-1}(\alpha)\bm{\sigma} ^+(\phi, \bm{\lambda}) =\mathbf{0},\\
& \bm{\lambda}\geq \mathbf{0},\,\, \mathbf{s}^-\geq \mathbf{0},\,\, \mathbf{s}^+\geq \mathbf{0},
\end{array}
\end{equation}
where $\phi =1+\beta $, and
\begin{align*}
\left( \sigma ^-_i\left( \bm{\lambda}\right) \right) ^2 =& \sum _{j,q=1}^n\lambda _j\lambda _q\mathrm{Cov}(\tilde{x}_{ij},\tilde{x}_{iq})-2\sum _{j=1}^n\lambda _j\mathrm{Cov}(\tilde{x}_{ij},\tilde{x}_{io})+\mathrm{Var}(\tilde{x}_{io}), \\
& i=1,\ldots ,m,\\
\left( \sigma ^+_r\left( \phi , \bm{\lambda}\right) \right) ^2 =& \sum _{j,q=1}^n\lambda _j\lambda _q\mathrm{Cov}(\tilde{y}_{rj},\tilde{y}_{rq})-2\phi \sum _{j=1}^n\lambda _j\mathrm{Cov}(\tilde{y}_{rj},\tilde{y}_{ro})+\phi ^2\mathrm{Var}(\tilde{y}_{ro}),\\
& r=1,\ldots ,s.
\end{align*}
Program \eqref{eq:ccr2d1o} is the chance constrained E-model corresponding to the output-oriented CCR model \eqref{eq:rad1b}. In fact, model \eqref{eq:ccr2d1o} (with unit weights) was constructed by \citet{Cooper2002}.

\begin{remark}
Both chance constrained radial models, the input-oriented \eqref{eq:ccr2d1i} and the output-oriented \eqref{eq:ccr2d1o}, were used by \citet{Lan1993} and \citet{Cooper2002}, respectively, to define concepts of stochastic efficiency. Nevertheless, taking into account Remark \ref{teo:2} and the fact that \eqref{eq:ccr2d1i} and \eqref{eq:ccr2d1o} are particular cases of \eqref{eq:ccdirD2d1} (that is equivalent to \eqref{eq:ccdirD21a}) with $\theta =1-\beta $ and $\phi =1+\beta $, we obtain that the concept of stochastic efficiency given in Definition \ref{def:2} coincides with those introduced by \citet{Lan1993} and \citet{Cooper2002}.
\end{remark}


Finally, if $D^-=I_m$ and $D^+=I_s$ (i.e. $d^-_1,\ldots ,d^-_m,d^+_1,\ldots ,d^+_s=1$), then we obtain a chance constrained E-model for computing a \textit{stochastic generalized Farrell measure} as the optimal value $\beta ^*$ of the following quadratic program:
\begin{equation*}
\label{eq:gfm}
\def\arraystretch{1.2}
\begin{array}[t]{rl}
\max \limits_{\beta ,\bm{\lambda}} & \beta \\
\text{s.t.} & (1-\beta )I_m\mathbf{x}_o-X\bm{\lambda}+\Phi ^{-1}(\alpha)\bm{\sigma} ^-(\beta ,\bm{\lambda}) \geq \mathbf{0},\\
& (1+\beta )I_s\mathbf{y}_{o}-Y\bm{\lambda}-\Phi ^{-1}(\alpha)\bm{\sigma} ^+(\beta ,\bm{\lambda}) \leq \mathbf{0},\\
& \bm{\lambda}\geq \mathbf{0},
\end{array}
\end{equation*}
where $\bm{\sigma} ^-(\beta ,\bm{\lambda})$ and $\bm{\sigma} ^+(\beta ,\bm{\lambda})$ are given by \eqref{eq:sigma1} and \eqref{eq:sigma2}, respectively.

\subsection{Deterministic directions}
\label{sec:ccddet}

In this section, we are going to apply the E-model approach to program \eqref{eq:dir1}, considering that inputs and outputs become stochastic but the directions $\mathbf{g}^-$, $\mathbf{g}^+$ remain deterministic. Hence, given $0<\alpha <1$, the corresponding chance constrained E-model can be written as
\begin{equation}
\label{eq:ccdirD21adet}
\def\arraystretch{1.8}
\begin{array}[t]{rl}
\max \limits_{\beta, \bm{\lambda},\mathbf{s}^-,\mathbf{s}^+} & \beta +\varepsilon \left( \mathbf{w}^-\mathbf{s}^-+\mathbf{w}^+\mathbf{s}^+\right) \\
\text{s.t.} & P\left\{ \left( \tilde{\mathbf{x}}_o-\beta \mathbf{g}^--\tilde{X}\bm{\lambda}-\mathbf{s}^-\right) _i\geq 0\right\} \geq 1-\alpha,\quad i=1,\ldots ,m, \\
& P\left\{ \left( \tilde{\mathbf{y}}_{o}+\beta \mathbf{g}^+-\tilde{Y}\bm{\lambda}+\mathbf{s}^+\right) _r\leq 0\right\} \geq 1-\alpha,\quad r=1,\ldots ,s, \\
& \bm{\lambda}\geq \mathbf{0},\,\, \mathbf{s}^-\geq \mathbf{0},\,\, \mathbf{s}^+\geq \mathbf{0}.
\end{array}
\end{equation}

\begin{remark}
The optimal objective value of \eqref{eq:ccdirD21adet} is $0$ if and only if the optimal objective value of \eqref{eq:ccdirD21a} is $0$, regardless of the directions. Hence, if we use program \eqref{eq:ccdirD21adet} instead of program \eqref{eq:ccdirD21a} in Definition \ref{def:2}, then we obtain an equivalent concept of stochastic efficiency. Analogously, the same applies to the concept of stochastic weak efficiency.
\end{remark}

If we assume that inputs and outputs are correlated random variables following multivariate normal distributions with means $E(\tilde{x}_{ij})=x_{ij}$ and $E(\tilde{y}_{rj})=y_{rj}$, we can deduce the deterministic equivalent to \eqref{eq:ccdirD21adet} following the methodology given in Section \ref{sec:ccdsto}:
\begin{equation}
\label{eq:ccdirD2d1det}
\def\arraystretch{1.2}
\begin{array}[t]{rl}
\max \limits_{\beta ,\bm{\lambda},\mathbf{s}^-,\mathbf{s}^+} & \beta +\varepsilon \left( \mathbf{w}^-\mathbf{s}^-+\mathbf{w}^+\mathbf{s}^+\right) \\
\text{s.t.} & \beta \mathbf{g}^-+X\bm{\lambda}+\mathbf{s}^--\Phi ^{-1}(\alpha)\bm{\sigma} ^-(\bm{\lambda}) = \mathbf{x}_o,\\
& -\beta \mathbf{g}^++Y\bm{\lambda}-\mathbf{s}^++\Phi ^{-1}(\alpha)\bm{\sigma} ^+(\bm{\lambda}) = \mathbf{y}_{o},\\
& \bm{\lambda}\geq \mathbf{0},\,\, \mathbf{s}^-\geq \mathbf{0},\,\, \mathbf{s}^+\geq \mathbf{0},
\end{array}
\end{equation}
where
\begin{align}
\label{eq:sigmadet1}
\left( \sigma ^-_i\left( \bm{\lambda}\right) \right) ^2 = & \sum _{j,q=1}^n\lambda _j\lambda _q\mathrm{Cov}(\tilde{x}_{ij},\tilde{x}_{iq})-2\sum _{j=1}^n\lambda _j\mathrm{Cov}(\tilde{x}_{ij},\tilde{x}_{io}) +
\mathrm{Var}(\tilde{x}_{io}),\\
\nonumber & i=1,\ldots ,m,\\
\label{eq:sigmadet2}
\left( \sigma ^+_r\left( \bm{\lambda}\right)\right) ^2 = & \sum _{j,q=1}^n\lambda _j\lambda _q\mathrm{Cov}(\tilde{y}_{rj},\tilde{y}_{rq})-2\sum _{j=1}^n\lambda _j\mathrm{Cov}(\tilde{y}_{rj},\tilde{y}_{ro})+
\mathrm{Var}(\tilde{y}_{ro}),\\
\nonumber & r=1,\ldots ,s.
\end{align}
Although model \eqref{eq:ccdirD2d1det} does not generalize chance constrained radial models due to the deterministic behaviour of its directions, it may be appropriate in certain cases, as we will show in Section \ref{sec:ex}.

\begin{remark}
\label{rem:ass}
Given the model \eqref{eq:ccdirD2d1} with stochastic directions of the form $\mathbf{g}^-=D^-\tilde{\mathbf{x}}_o$ and $\mathbf{g}^+=D^+\tilde{\mathbf{y}}_o$, 
we can consider a model \eqref{eq:ccdirD2d1det} with deterministic directions given by the means of the stochastic directions, i.e. $\mathbf{g}^-=D^-\mathbf{x}_o$ and $\mathbf{g}^+=D^+\mathbf{y}_o$. We can write this model as:
\begin{equation}
\label{eq:ccdirD2d1ass}
\def\arraystretch{1.2}
\begin{array}[t]{rl}
\max \limits_{\beta ,\bm{\lambda},\mathbf{s}^-,\mathbf{s}^+} & \beta +\varepsilon \left( \mathbf{w}^-\mathbf{s}^-+\mathbf{w}^+\mathbf{s}^+\right) \\
\text{s.t.} & \Theta ^-(\beta )\mathbf{x}_o-X\bm{\lambda}-\mathbf{s}^-+\Phi ^{-1}(\alpha)\bm{\sigma} ^-(\bm{\lambda}) = \mathbf{0},\\
& \Theta ^+(\beta )\mathbf{y}_{o}-Y\bm{\lambda}+\mathbf{s}^+-\Phi ^{-1}(\alpha)\bm{\sigma} ^+(\bm{\lambda}) = \mathbf{0},\\
& \bm{\lambda}\geq \mathbf{0},\,\, \mathbf{s}^-\geq \mathbf{0},\,\, \mathbf{s}^+\geq \mathbf{0},
\end{array}
\end{equation}
where $\bm{\sigma} ^+(\bm{\lambda})$ and $\bm{\sigma} ^-(\bm{\lambda})$ are given by \eqref{eq:sigmadet1} and \eqref{eq:sigmadet2}, respectively. In this case, we say that both models, \eqref{eq:ccdirD2d1} and \eqref{eq:ccdirD2d1ass}, are \emph{associated}. We have that associated models are equivalent in the deterministic case, i.e. if the variances of the stochastic variables are $0$.

Conversely, given the model \eqref{eq:ccdirD2d1det} with deterministic directions $\mathbf{g}^-$ and $\mathbf{g}^+$, we can always find an associated model \eqref{eq:ccdirD2d1} with stochastic directions of the form $\mathbf{g}^-=D^-\tilde{\mathbf{x}}_o$ and $\mathbf{g}^+=D^+\tilde{\mathbf{y}}_o$, where $D^-=\mathrm{diag}(d^-_1,\ldots ,d^-_m)$ and $D^+=\mathrm{diag}(d^+_1,\ldots ,d^+_s)$ with $d_i^-=g_i^-/x_{io}$, $d_r^+=g_r^+/y_{ro}$.
\end{remark}

\section{Examples}
\label{sec:ex}

We are going to consider the example given in \citet{Charnes1981, Lan1993}, in which the ``Program Follow Through'' experiment is analyzed in public school education. We have to note that this example is given for illustrative purposes only. Some other meaningful examples could be given by considering factories instead of schools, workers instead of students, etc.
The original data set consists on 49 school sites enrolled in the experiment, with 5 inputs (education, level of mother, parent occupation, parental visit index, counseling index, and number of teachers) and 3 outputs (total reading scores, total math scores, and total Coopersmith\footnote{An index of a child's self-esteem.} scores). The average values of input and output variables are reported for each school site. Table \ref{tab:datos} shows these values for the first $10$ school sites.

\begin{table}[]
\begin{center}
\caption{Average values from the first $10$ school sites in the Program Follow Through.}
\label{tab:datos}
\begin{tabular}{lcc}
 & Inputs & Outputs \\
Site 1 & $(86.13,\, 16.24,\, 48.21,\, 49.69,\, 9)$ & $(54.53,\, 58.98,\, 38.16)$ \\
Site 2 & $(29.26,\, 10.24,\, 41.96,\, 40.65,\, 5)$ & $(24.69,\, 33.89,\, 26.02)$ \\
Site 3 & $(43.12,\, 11.31,\, 38.19,\, 35.03,\, 9)$ & $(36.41,\, 40.62,\, 28.51)$ \\
Site 4 & $(24.96,\, \,\,\,6.14,\, 24.81,\, 25.15,\, 7)$ & $(14.94,\, 17.58,\, 16.19)$ \\
Site 5 & $(11.62,\, \,\,\,2.21,\, \,\,\,6.85,\, \,\,\,6.37,\, 4)$ & $(\,\,\,7.81,\, \,\,\,6.94,\, \,\,\,5.37)$ \\
Site 6 & $(11.88,\, \,\,\,4.97,\, 18.73,\, 18.04,\, 4)$ & $(12.59,\, 16.85,\, 12.84)$ \\
Site 7 & $(32.64,\, \,\,\,6.88,\, 28.10,\, 25.45,\, 7)$ & $(17.06,\, 16.99,\, 17.82)$ \\
Site 8 & $(20.79,\, 12.97,\, 54.85,\, 52.07,\, 8)$ & $(20.19,\, 30.64,\, 33.16)$ \\
Site 9 & $(34.40,\, 11.04,\, 38.16,\, 42.40,\, 8)$ & $(26.13,\, 29.80,\, 26.29)$ \\
Site 10 & $(61.74,\, 14.50,\, 49.09,\, 42.92,\, 9)$ & $(46.42,\, 51.59,\, 35.20)$      
\end{tabular}
\end{center}
\end{table}

In order to simplify the problem, we are going to assume that inputs are deterministic and outputs are stochastic. Moreover, we assume that all observed outputs coincide with their mathematical expectations, all outputs are stochastically independent, and the within-school variability of each output (as measured by the variance) is the same, $c^2$, for all outputs and at all school sites. See \citet{Lan1993} for a detailed discussion on these assumptions.

In this example, we are interested in knowing how much each school site's output scores must improve to become $\alpha$-stochastically efficient or weakly efficient and therefore, only the optimal score $\beta ^*$ of  an output-oriented model (i.e. with $\mathbf{g}^-=\mathbf{0}$ or $D^-=0$) will be necessary. The rationale for considering weak efficiency as a goal (jointly with efficiency) lies in the fact that any small improvement in some output that has no slacks, results in the DMU being $\alpha$-stochastically efficient. Note that this is always possible since, in output-oriented models, there is always some output without slacks.

If we take $D^-=0$ and $D^+=I_3$ (i.e. $d_1^+=d_2^+=d_3^+=1$) in program \eqref{eq:ccdirD2d1}, then the chance constrained output-oriented CCR model \eqref{eq:ccr2d1o} is applied. However, this choice treats all evaluated subjects (reading, math and Coopersmith) in the same way and hence, it may not reflect the specific characteristics of improvement of each subject. On the other hand, if for example we set $d_1^+=0.1$, $d_2^+=0.05$ and $d_3^+=0.01$, then we suppose that students' ability to improve the reading score by $10\%$ is the same as their ability to improve the math score by $5\%$, and the same as their ability to improve the Coopersmith score by $1\%$.

If we look at improvement capabilities in absolute terms, then the directions become deterministic and we must apply model \eqref{eq:ccdirD2d1det} with $\mathbf{g}^-=\mathbf{0}$. For example, if we set $g_1^+=5$, $g_2^+=4$ and $g_3^+=1$, then we suppose that students' ability to improve the reading score by $5$ absolute points is the same as their ability to improve the math score by $4$ absolute points, and the same as their ability to improve the Coopersmith score by $1$ absolute point.

Note that, in order to compare optimal scores $\beta ^*$ from models with different directions, we have to take into account that those directions should represent the same ``amount of effort'' on the part of the students of the corresponding evaluated school site. For example, if the students in the Site 1 can improve their reading scores by $10\%$, their math scores by $5\%$ and their Coopersmith scores by $1\%$ using $1$ ``unit of effort'', and the students in the Site 2 can improve their reading scores by $5$ points, their math scores by $4$ points and their Coopersmith scores by $1$ point also using $1$ ``unit of effort'', then we can compare the $\beta ^*$ scores of both sites if we use model \eqref{eq:ccdirD2d1} with $D^-=0$, $d_1^+=0.1$, $d_2^+=0.05$, $d_3^+=0.01$ for evaluating Site 1, and model \eqref{eq:ccdirD2d1det} with $\mathbf{g}^-=\mathbf{0}$, $g_1^+=5$, $g_2^+=4$, $g_3^+=1$ for evaluating Site 2. In this case, $\beta ^*$ is interpreted as the ``total amount of effort'' of the corresponding school site for becoming  $\alpha$-stochastically efficient or weakly efficient.

\begin{table}[]
\begin{center}
\caption{Results of the optimal $\beta ^*$ scores of chance constrained directional models with stochastic directions ($D^-=0$, $d_1^+=1$, $d_2^+=1$, $d_3^+=1$) and their associated models with deterministic directions ($\mathbf{g}^-=\mathbf{0}$, $g_1^+=y_{1o}$, $g_2^+=y_{2o}$, $g_3^+=y_{3o}$), applied to the first $10$ school sites using different values $c^2$ of the variance of the outputs and $\alpha = 0.05$.}

$ $

\label{tab:1}
\begin{tabular}{lcccccccccc}
\hline
\multicolumn{11}{c}{Stochastic directions: $D^-=0,\quad d_1^+=1,\quad d_2^+=1,\quad d_3^+=1$}  \\ \hline
 Site & 1 & 2 & 3 & 4 & 5 & 6 & 7 & 8 & 9 & 10 \\ \hline
 $c=0$ & 0 & 0.109 & 0.012 & 0.108 & 0 & 0.103 & 0.121 & 0.093 & 0.148 & 0 \\
 $c=0.5$ & 0 & 0.071 & 0 & 0.042 & 0 & 0.031 & 0.061 & 0.063 & 0.095 & 0 \\ 
 $c=1$ & 0 & 0.036 & 0 & 0 & 0 & 0 & 0.006 & 0.026 & 0.053 & 0 \\ \hline
 \hline
 \multicolumn{11}{c}{Deterministic directions: $\mathbf{g}^-=\mathbf{0},\quad g_1^+=y_{1o},\quad g_2^+=y_{2o},\quad g_3^+=y_{3o}$}  \\ \hline
 Site & 1 & 2 & 3 & 4 & 5 & 6 & 7 & 8 & 9 & 10 \\ \hline
 $c=0$ & 0 & 0.109 & 0.012 & 0.108 & 0 & 0.103 & 0.121 & 0.093 & 0.148 & 0 \\
 $c=0.5$ & 0 & 0.073 & 0 & 0.044 & 0 & 0.033 & 0.063 & 0.065 & 0.098 & 0 \\
 $c=1$ & 0 & 0.038 & 0 & 0 & 0 & 0 & 0.007 & 0.033 & 0.055 & 0 \\ \hline   
\end{tabular}
\end{center}
\end{table}

Taking all this into account, we have evaluated the first $10$ school sites with respect to the whole sample, for $\alpha =0.05$ and using different values $c^2$ of the variance of the outputs, obtaining the deterministic case by taking $c=0$. Tables \ref{tab:1}, \ref{tab:2} and \ref{tab:3} show the results of the optimal $\beta ^*$ scores of chance constrained directional models with stochastic directions (columns $1-3$) and their associated models with deterministic directions (columns $4-6$). Specifically, Table \ref{tab:1} shows the results of applying the output-oriented chance constrained CCR model \eqref{eq:ccr2d1o}, that does not reflect the specific characteristics of improvement of each subject. Analogously, Table \ref{tab:2} shows the results of applying the chance constrained directional model \eqref{eq:ccdirD2d1} with stochastic directions determined by orientation parameters $D^-=0$, $d_1^+=0.1$, $d_2^+=0.05$ and $d_3^+=0.01$. In this case, we are assuming that the students can improve the reading score by $10\%$ the math score by $5\%$ and the Coopersmith score by $1\%$ using $1$ ``unit of effort'' in all school sites. Finally, Table \ref{tab:3} shows the results of applying the chance constrained directional model \eqref{eq:ccdirD2d1det} with deterministic directions given by $\mathbf{g}^-=\mathbf{0}$, $g_1^+=5$, $g_2^+=4$ and $g_3^+=1$, assuming that the students can improve the reading score by $5$ points, the math score by $4$ points and the Coopersmith score by $1$ point using $1$ ``unit of effort'' in all school sites.

\begin{table}[]
\begin{center}
\caption{Results of the optimal $\beta ^*$ scores of chance constrained directional models with stochastic directions ($D^-=0$, $d_1^+=0.1$, $d_2^+=0.05$, $d_3^+=0.01$) and their associated models with deterministic directions ($\mathbf{g}^-=\mathbf{0}$, $g_1^+=d_1^+y_{1o}$, $g_2^+=d_2^+y_{2o}$, $g_3^+=d_3^+y_{3o}$), applied to the first $10$ school sites using different values $c^2$ of the variance of the outputs and $\alpha = 0.05$.}

$ $ 

\label{tab:2}
\begin{tabular}{lcccccccccc}
\hline
\multicolumn{11}{c}{Stochastic directions: $D^-=0,\quad d_1^+=0.1,\quad d_2^+=0.05,\quad d_3^+=0.01$}  \\ \hline
 Site & 1 & 2 & 3 & 4 & 5 & 6 & 7 & 8 & 9 & 10 \\ \hline
 $c=0$ & 0 & 5.041 & 0.388 & 4.988 & 0 & 3.380 & 5.468 & 8.218 & 5.303 & 0 \\
 $c=0.5$ & 0 & 3.601 & 0 & 2.117 & 0 & 1.664 & 2.876 & 6.301 & 4.481 & 0 \\
 $c=1$ & 0 & 2.296 & 0 & 0 & 0 & 0 & 0.374 & 3.409 & 3.573 & 0 \\ \hline
 \hline
 \multicolumn{11}{c}{Deterministic directions: $\mathbf{g}^-=\mathbf{0},\,\,\, g_1^+=d_1^+y_{1o},\,\,\, g_2^+=d_2^+y_{2o},\,\,\, g_3^+=d_3^+y_{3o}$}  \\ \hline
 Site & 1 & 2 & 3 & 4 & 5 & 6 & 7 & 8 & 9 & 10 \\ \hline
 $c=0$ & 0 & 5.041 & 0.388 & 4.988 & 0 & 3.380 & 5.468 & 8.218 & 5.303 & 0 \\
 $c=0.5$ & 0 & 3.707 & 0 & 2.216 & 0 & 1.768 & 2.994 & 6.437 & 4.592 & 0 \\
 $c=1$ & 0 & 2.426 & 0 & 0 & 0 & 0 & 0.404 & 3.555 & 3.752 & 0 \\ \hline  
\end{tabular}
\end{center}
\end{table}

\begin{table}[]
\begin{center}
\caption{Results of the optimal $\beta ^*$ scores of chance constrained directional models with  deterministic directions ($\mathbf{g}^-=\mathbf{0}$, $g_1^+=5$, $g_2^+=4$, $g_3^+=1$) and their associated models with stochastic directions ($D^-=0$, $d_1^+=5/y_{1o}$, $d_2^+=4/y_{2o}$, $d_3^+=1/y_{3o}$), applied to the first $10$ school sites using different values $c^2$ of the variance of the outputs and $\alpha = 0.05$.}

$ $ 

\label{tab:3}
\begin{tabular}{lcccccccccc}
\hline
\multicolumn{11}{c}{Stochastic directions: $D^-=0,\quad d_1^+=5/y_{1o},\quad d_2^+=4/y_{2o},\quad d_3^+=1/y_{3o}$}  \\ \hline
 Site & 1 & 2 & 3 & 4 & 5 & 6 & 7 & 8 & 9 & 10 \\ \hline
 $c=0$ & 0 & 1.982 & 0.211 & 1.137 & 0 & 0.754 & 1.412 & 3.090 & 2.561 & 0 \\
 $c=0.5$ & 0 & 1.415 & 0 & 0.466 & 0 & 0.338 & 0.729 & 2.089 & 2.100 & 0 \\
 $c=1$ & 0 & 0.819 & 0 & 0 & 0 & 0 & 0.080 & 1.130 & 1.413 & 0 \\ \hline
 \hline
 \multicolumn{11}{c}{Deterministic directions: $\mathbf{g}^-=\mathbf{0},\quad g_1^+=5,\quad g_2^+=4,\quad g_3^+=1$}  \\ \hline
 Site & 1 & 2 & 3 & 4 & 5 & 6 & 7 & 8 & 9 & 10 \\ \hline
 $c=0$ & 0 & 1.982 & 0.211 & 1.137 & 0 & 0.754 & 1.412 & 3.090 & 2.561 & 0 \\
 $c=0.5$ & 0 & 1.457 & 0 & 0.487 & 0 & 0.359 & 0.755 & 2.134 & 2.152 & 0 \\
 $c=1$ & 0 & 0.864 & 0 & 0 & 0 & 0 & 0.093 & 1.179 & 1.483 & 0 \\ \hline  
\end{tabular}
\end{center}
\end{table}

Regarding the discussion of the results, first of all, associated models with stochastic and deterministic directions coincide in the deterministic case ($c=0$), as stated in Remark \ref{rem:ass}. Second, the greater the stochastic variability of outputs (the greater the coefficient $c$), the smaller the optimal score $\beta ^*$ is, and hence, school sites are less inefficient. The reason is that in the chance constrained formulation of DEA, the ``hard'' efficient frontier of deterministic DEA is replaced with a ``soft'' frontier that moves successively closer to any given observation \citep{Lan1993}. This effect is also seen among associated models, where models with stochastic directions produce better scores $\beta ^*$ than the corresponding associated models with deterministic directions. Third, the relative rankings of inefficient school sites differ in the chance constrained analysis as compared to deterministic DEA. For example, in Tables \ref{tab:1} and \ref{tab:2}, Site 7 is more inefficient than Site 2 for $c=0$, but this situation is reversed for $c=0.5$ and $c=1$. Moreover, in Table \ref{tab:3}, Site 8 is more inefficient than Site 9 for $c=0$, but the opposite occurs for $c=0.5$ and $c=1$. Finally, different directions (representing different improvement strategies) also produce different rankings, as expected. For example, in Table \ref{tab:1}, Site 2 is more inefficient than Site 8 for all values of $c$, just the opposite of what happens in Tables \ref{tab:2} and \ref{tab:3}.  

We have used R 4.2.0 \citep{R22} for computations. Specifically, we have used the optiSolve package \citep{optiSolve} for solving quadratic programs, and the deaR package \citep{deaR22} for linear models, in order to check the deterministic case $c=0$.

\section{Concluding remarks}
\label{sec:con}

We have constructed different versions of stochastic directional models following the chance constrained methodology and generalizing the existing radial models developed by \citet{Cooper2002}. The advantage of directional models over radial models is the great versatility provided by the choice of the direction, allowing the user to customize an improvement strategy that can be adapted to market conditions. Moreover, chance constrained models allow a closer approximation to a reality in which measurements present uncertainty, providing an alternative to fuzzy models.

Nevertheless, chance constrained DEA also has drawbacks. The major operational disadvantage is that programs are quadratic, as opposed to deterministic models whose programs are usually linear. However, quadratic programs are relatively fast to solve compared with general non linear programs. On the other hand, 
no formal statistical model with a sampling process is specified, unlike what is done by \citet{Banker1993}. The main consequence is the difficulty in defining an appropriate and meaningful measure of inefficiency \citep{Olesen2016}. Finally, according to Remark \ref{rem:multi}, a more appropriate chance constrained methodology would be to require that all the constraints are met at once with a certain probability.
But, in order to find the corresponding deterministic equivalent program, we should consider independent variables (thus significantly reducing the applicability), or use multivariate probability distributions, for which the inverse distribution functions are not uniquely defined. However, despite the difficulty of working with joint chance constrained models, they could lead to future studies and more sophisticated models in the chance constrained DEA framework.

\bibliography{refs}

\begin{thebibliography}{31}
\providecommand{\natexlab}[1]{#1}
\providecommand{\url}[1]{\texttt{#1}}
\expandafter\ifx\csname urlstyle\endcsname\relax
  \providecommand{\doi}[1]{doi: #1}\else
  \providecommand{\doi}{doi: \begingroup \urlstyle{rm}\Url}\fi

\bibitem[Amirteimoori et~al.(2022)Amirteimoori, Sahoo, Charles, and
  Mehdizadeh]{Amir2022}
A.~Amirteimoori, B.~K. Sahoo, V.~Charles, and S.~Mehdizadeh.
\newblock \emph{Stochastic Benchmarking}.
\newblock Springer Cham, 1st edition, 2022.
\newblock ISBN 978-3-030-89871-7.
\newblock \doi{10.1007/978-3-030-89869-4}.

\bibitem[Aparicio et~al.(2007)Aparicio, Ruiz, and Sirvent]{Apa2007}
J.~Aparicio, J.~L. Ruiz, and I.~Sirvent.
\newblock Closest targets and minimum distance to the {Pareto}-efficient
  frontier in {DEA}.
\newblock \emph{Journal of Productivity Analysis}, 28\penalty0 (3):\penalty0
  209--218, 2007.
\newblock \doi{10.1007/s11123-007-0039-5}.

\bibitem[Banker(1993)]{Banker1993}
R.~D. Banker.
\newblock Maximum likelihood, consistency and data envelopment analysis: A
  statistical foundation.
\newblock \emph{Management Science}, 39\penalty0 (10):\penalty0 1265--1273,
  1993.
\newblock \doi{10.1287/mnsc.39.10.1265}.

\bibitem[Banker et~al.(1984)Banker, Charnes, and Cooper]{Banker1984}
R.~D. Banker, A.~Charnes, and W.~W. Cooper.
\newblock Some models for estimating technical and scale inefficiencies in data
  envelopment analysis.
\newblock \emph{Management Science}, 30\penalty0 (9):\penalty0 1078--1092,
  1984.
\newblock \doi{10.1287/mnsc.30.9.1078}.

\bibitem[Briec(1997)]{Briec1997}
W.~Briec.
\newblock A graph-type extension of {Farrell} technical efficiency measure.
\newblock \emph{Journal of Productivity Analysis}, 8\penalty0 (1):\penalty0
  95--110, 1997.
\newblock \doi{10.1023/A:1007728515733}.

\bibitem[Chambers et~al.(1996)Chambers, Chung, and Färe]{Chambers1996}
R.~G. Chambers, Y.~Chung, and R.~Färe.
\newblock Benefit and distance functions.
\newblock \emph{Journal of Economic Theory}, 70\penalty0 (2):\penalty0
  407--419, 1996.
\newblock \doi{10.1006/jeth.1996.0096}.

\bibitem[Chambers et~al.(1998)Chambers, Chung, and F{\"{a}}re]{Chambers1998}
R.~G. Chambers, Y.~Chung, and R.~F{\"{a}}re.
\newblock Profit, directional distance functions, and {Nerlovian} efficiency.
\newblock \emph{Journal of Optimization Theory and Applications}, 98\penalty0
  (2):\penalty0 351--364, 1998.
\newblock \doi{10.1023/A:1022637501082}.

\bibitem[Charnes et~al.(1978)Charnes, Cooper, and Rhodes]{Charnes1978}
A.~Charnes, W.~Cooper, and E.~Rhodes.
\newblock Measuring the efficiency of decision making units.
\newblock \emph{European Journal of Operational Research}, 2\penalty0
  (6):\penalty0 429--444, 1978.
\newblock \doi{10.1016/0377-2217(78)90138-8}.

\bibitem[Charnes et~al.(1979)Charnes, Cooper, and Rhodes]{Charnes1979}
A.~Charnes, W.~Cooper, and E.~Rhodes.
\newblock Short communication: Measuring the efficiency of decision making
  units.
\newblock \emph{European Journal of Operational Research}, 3\penalty0
  (4):\penalty0 339, 1979.
\newblock \doi{10.1016/0377-2217(79)90229-7}.

\bibitem[Charnes et~al.(1981)Charnes, Cooper, and Rhodes]{Charnes1981}
A.~Charnes, W.~W. Cooper, and E.~Rhodes.
\newblock Evaluating program and managerial efficiency: an application of data
  envelopment analysis to {Program Follow Through}.
\newblock \emph{Management Science}, 27\penalty0 (6):\penalty0 668--697, 1981.
\newblock \doi{10.1287/mnsc.27.6.668}.

\bibitem[Charnes et~al.(1985)Charnes, Cooper, Golany, Seiford, and
  Stutz]{Charnes1985}
A.~Charnes, W.~Cooper, B.~Golany, L.~Seiford, and J.~Stutz.
\newblock Foundations of data envelopment analysis for {Pareto-Koopmans}
  efficient empirical production functions.
\newblock \emph{Journal of Econometrics}, 30\penalty0 (1-2):\penalty0 91--107,
  1985.
\newblock \doi{10.1016/0304-4076(85)90133-2}.

\bibitem[Coll-Serrano et~al.(2022)Coll-Serrano, Bol\'os, and
  Ben\'{\i}tez]{deaR22}
V.~Coll-Serrano, V.~J. Bol\'os, and R.~Ben\'{\i}tez.
\newblock \emph{{deaR}: Conventional and Fuzzy Data Envelopment Analysis},
  2022.
\newblock URL \url{https://CRAN.R-project.org/package=deaR}.
\newblock R package version 1.3.3.

\bibitem[Cooper et~al.(1996)Cooper, Huang, and Li]{Cooper1996}
W.~W. Cooper, Z.~M. Huang, and S.~X. Li.
\newblock Satisficing {DEA} models under chance constraints.
\newblock \emph{Annals of Operations Research}, 66:\penalty0 279--295, 1996.
\newblock \doi{10.1007/BF02187302}.

\bibitem[Cooper et~al.(1998)Cooper, Huang, Lelas, Li, and Olesen]{Cooper1998}
W.~W. Cooper, Z.~Huang, V.~Lelas, S.~X. Li, and O.~B. Olesen.
\newblock Chance constrained programming formulations for stochastic
  characterizations of efficiency and dominance in {DEA}.
\newblock \emph{Journal of Productivity Analysis}, 9\penalty0 (1):\penalty0
  53--79, 1998.
\newblock \doi{10.1023/A:1018320430249}.

\bibitem[Cooper et~al.(2002)Cooper, Deng, Huang, and Li]{Cooper2002}
W.~W. Cooper, H.~Deng, Z.~Huang, and S.~X. Li.
\newblock Chance constrained programming approaches to technical efficiencies
  and inefficiencies in stochastic data envelopment analysis.
\newblock \emph{Journal of the Operational Research Society}, 53\penalty0
  (12):\penalty0 1347--1356, 2002.
\newblock \doi{10.1057/palgrave.jors.2601433}.

\bibitem[Cooper et~al.(2007)Cooper, Seiford, and Tone]{Cooper2007}
W.~W. Cooper, L.~M. Seiford, and K.~Tone.
\newblock \emph{Data Envelopment Analysis. A Comprehensive Text with Models,
  Applications, References and {DEA}-Solver Software}.
\newblock Springer, 2nd edition, 2007.
\newblock ISBN 9780387452814.
\newblock \doi{10.1007/978-0-387-45283-8}.

\bibitem[F{\"a}re and {Knox Lovell}(1978)]{Fare1978}
R.~F{\"a}re and C.~A. {Knox Lovell}.
\newblock Measuring the technical efficiency of production.
\newblock \emph{Journal of Economic Theory}, 19\penalty0 (1):\penalty0
  150--162, 1978.
\newblock \doi{10.1016/0022-0531(78)90060-1}.

\bibitem[Guo and Tanaka(2001)]{Guo2001}
P.~Guo and H.~Tanaka.
\newblock Fuzzy {DEA}: A perceptual evaluation method.
\newblock \emph{Fuzzy Sets and Systems}, 119\penalty0 (1):\penalty0 149--160,
  2001.
\newblock \doi{10.1016/S0165-0114(99)00106-2}.

\bibitem[Kao and Liu(2000)]{Kao2000}
C.~Kao and S.-T. Liu.
\newblock Fuzzy efficiency measures in data envelopment analysis.
\newblock \emph{Fuzzy Sets and Systems}, 113\penalty0 (3):\penalty0 427--437,
  2000.
\newblock \doi{10.1016/S0165-0114(98)00137-7}.

\bibitem[Land et~al.(1993)Land, Knox~Lovell, and Thore]{Lan1993}
K.~C. Land, C.~A. Knox~Lovell, and S.~Thore.
\newblock Chance-constrained data envelopment analysis.
\newblock \emph{Managerial and Decision Economics}, 14\penalty0 (6):\penalty0
  541--554, 1993.
\newblock https://www.jstor.org/stable/2487873.

\bibitem[Le{\'{o}}n et~al.(2003)Le{\'{o}}n, Liern, Ruiz, and Sirvent]{Leon2003}
T.~Le{\'{o}}n, V.~Liern, J.~L. Ruiz, and I.~Sirvent.
\newblock A fuzzy mathematical programming approach to the assessment of
  efficiency with {DEA} models.
\newblock \emph{Fuzzy Sets and Systems}, 139\penalty0 (2):\penalty0 407--419,
  2003.
\newblock \doi{10.1016/S0165-0114(02)00608-5}.

\bibitem[Luenberger(1992)]{Luenberger1992}
D.~G. Luenberger.
\newblock New optimality principles for economic efficiency and equilibrium.
\newblock \emph{Journal of Optimization Theory and Applications}, 75\penalty0
  (2):\penalty0 221--264, 1992.
\newblock \doi{10.1007/BF00941466}.

\bibitem[Olesen and Petersen(1995)]{Olesen1995}
O.~B. Olesen and N.~C. Petersen.
\newblock Chance constrained efficiency evaluation.
\newblock \emph{Management Science}, 41\penalty0 (3):\penalty0 442--457, 1995.
\newblock \doi{10.1287/mnsc.41.3.442}.

\bibitem[Olesen and Petersen(2016)]{Olesen2016}
O.~B. Olesen and N.~C. Petersen.
\newblock Stochastic data envelopment analysis—a review.
\newblock \emph{European Journal of Operational Research}, 251\penalty0
  (1):\penalty0 2--21, 2016.
\newblock \doi{10.1016/j.ejor.2015.07.058}.

\bibitem[Pastor et~al.(1999)Pastor, Ruiz, and Sirvent]{Pastor1999}
J.~Pastor, J.~Ruiz, and I.~Sirvent.
\newblock An enhanced {DEA} {Russell} graph efficiency measure.
\newblock \emph{European Journal of Operational Research}, 115\penalty0
  (3):\penalty0 596--607, 1999.
\newblock ISSN 0377-2217.
\newblock \doi{10.1016/S0377-2217(98)00098-8}.

\bibitem[{R Core Team}(2022)]{R22}
{R Core Team}.
\newblock \emph{R: A Language and Environment for Statistical Computing}.
\newblock R Foundation for Statistical Computing, Vienna, Austria, 2022.
\newblock URL \url{https://www.R-project.org/}.

\bibitem[Shiraz et~al.(2021)Shiraz, Tavana, and Fukuyama]{Shiraz2021}
R.~K. Shiraz, M.~Tavana, and H.~Fukuyama.
\newblock A joint chance-constrained data envelopment analysis model with
  random output data.
\newblock \emph{Operational Research}, 21:\penalty0 1255--1277, 2021.
\newblock \doi{10.1007/s12351-019-00478-0}.

\bibitem[Tavana et~al.(2013)Tavana, Shiraz, Hatami-Marbini, Agrell, and
  Paryab]{Tavana2013}
M.~Tavana, R.~K. Shiraz, A.~Hatami-Marbini, P.~Agrell, and K.~Paryab.
\newblock Chance-constrained dea models with random fuzzy inputs and outputs.
\newblock \emph{Knowledge-Based Systems}, 52:\penalty0 32--52, 2013.
\newblock \doi{10.1016/j.knosys.2013.05.014}.

\bibitem[Tone(2001)]{Tone2001}
K.~Tone.
\newblock A slacks-based measure of efficiency in data envelopment analysis.
\newblock \emph{European Journal of Operational Research}, 130\penalty0
  (3):\penalty0 498--509, 2001.
\newblock \doi{10.1016/S0377-2217(99)00407-5}.

\bibitem[Tone(2010)]{Tone2010}
K.~Tone.
\newblock Variations on the theme of slacks-based measure of efficiency in
  {DEA}.
\newblock \emph{European Journal of Operational Research}, 200\penalty0
  (3):\penalty0 901--907, 2010.
\newblock \doi{10.1016/j.ejor.2009.01.027}.

\bibitem[Wellmann(2021)]{optiSolve}
R.~Wellmann.
\newblock \emph{{optiSolve}: Linear, Quadratic, and Rational Optimization},
  2021.
\newblock URL \url{https://CRAN.R-project.org/package=optiSolve}.
\newblock R package version 1.0.

\end{thebibliography}
\bibliographystyle{abbrvnat}

\end{document}